\documentclass{amsart}%
\usepackage{amssymb}
\usepackage{amsfonts}
\usepackage{amsmath}
\usepackage{graphicx}%
\newtheorem{theorem}{Theorem}
\theoremstyle{plain}

\newtheorem{lemma}{Lemma}

\newtheorem{remark}{Remark}

\numberwithin{equation}{section}
\begin{document}
\title[ ]{ A Sharp Inequality Relating Yamabe Invariants on Asymptotically
Poincare-Einstein Manifolds with a Ricci Curvature Lower Bound}
\author{Xiaodong Wang}
\address{Department of Mathematics, Michigan State University, East Lansing, MI 48824}
\email{xwang@msu.edu}
\author{Zhixin Wang}
\address{Department of Mathematics, Michigan State University, East Lansing, MI 48824}
\email{wangz117@msu.edu}

\begin{abstract}
Let $\left(  X^{n},g_{+}\right)  $ be a conformally compact manifold with
$Ric\geq-\left(  n-1\right)  $. If $g_{+}$ is asymptotically
Poincare-Einstein, we establish a sharp inequality relating the type II Yamabe
invariant of $\overline{X}$ and the Yamabe invariant of its conformal infinity.

\end{abstract}
\maketitle

\section{Introduction}

The Yamabe problem for closed Riemannian manifolds was completely solved by
Aubin and Schoen (cf. \cite{A, SY} for complete exposition). For compact
Riemannian manifolds with boundary, there are two types of Yamabe problems and
neither has bee completely solved. Let $\left(  M^{n},g\right)  $ be a compact
Riemannian manifold $\left(  M^{n},g\right)  $ with nonempty boundary
$\Sigma=\partial M$. The functional%
\[
E_{g}\left(  u\right)  =\int_{M}\left(  \frac{4\left(  n-1\right)  }%
{n-2}\left\vert \nabla u\right\vert ^{2}+Ru^{2}\right)  dv_{g}+2\int_{\Sigma
}Hu^{2}d\sigma_{g},
\]
where $R$ is the scalar curvature and $H$ is the mean curvature of the
boundary, has the important property of being conformally invariant: if
$\widetilde{g}=\phi^{4/\left(  n-2\right)  }g$ is another metric, then
$E_{\widetilde{g}}\left(  u\right)  =E_{g}\left(  u\phi\right)  $. The
functional can be written as%
\[
E_{g}\left(  u\right)  =\int_{M}uL_{g}udv_{g}+2\int_{\Sigma}\left(
\frac{2\left(  n-1\right)  }{n-2}\frac{\partial u}{\partial\nu}+Hu\right)
ud\sigma_{g},
\]
where $L_{g}u=-\frac{4\left(  n-1\right)  }{n-2}\Delta_{g}u+Ru$ is the
conformal Laplacian. The type I Yamabe invariant is defined as
\[
Y\left(  M,\left[  g\right]  \right)  =\inf_{u\in H^{1}\left(  M\right)
\backslash\left\{  0\right\}  }\frac{E_{g}\left(  u\right)  }{\left(  \int
_{M}\left\vert u\right\vert ^{2n/\left(  n-2\right)  }dv_{g}\right)  ^{\left(
n-2\right)  /n}}.
\]
The type I Yamabe problem is whether the infimum is always achieved. It is
proved that $Y\left(  M,\left[  g\right]  \right)  \leq Y\left(
\mathbb{S}_{+}^{n}\right)  $ and moreover the infimum is achieved if the
inequality is strict. The strategy to solve the type I Yamabe problem is to
show that the strict inequality $Y\left(  M,\left[  g\right]  \right)
<Y\left(  \mathbb{S}_{+}^{n}\right)  $ is always true unless $\left(
M,\left[  g\right]  \right)  $ is conformal diffeomorphic to $\mathbb{S}%
_{+}^{n}$. It has been confirmed in many cases (see \cite{E1} and \cite{BC}),
but some exceptional cases remain open.

The type II Yamabe invariant is defined as
\[
Q\left(  M,\Sigma,\left[  g\right]  \right)  =\inf_{u\in H^{1}\left(
M\right)  \backslash\left\{  0\right\}  }\frac{E_{g}\left(  u\right)
}{\left(  \int_{\Sigma}\left\vert u\right\vert ^{2\left(  n-1\right)  /\left(
n-2\right)  }d\sigma_{g}\right)  ^{\left(  n-2\right)  /\left(  n-1\right)  }%
}.
\]
It should be noted that $Q\left(  M,\Sigma,\left[  g\right]  \right)  $ can be
$-\infty$. If $Q\left(  M,\Sigma,\left[  g\right]  \right)  $ $>-\infty$ and
the infimum is achieved, then a minimizer $u$ properly scaled is smooth and
positive and the metric $u^{4/\left(  n-2\right)  }g$ then has zero scalar
curvature on $M$ and constant mean curvature on $\Sigma$. The type II Yamabe
problem is whether the infimum is achieved when $Q\left(  M,\Sigma,\left[
g\right]  \right)  $ $>-\infty$. Parallel to the type I Yamabe problem, it is
proved that $Q\left(  M,\Sigma,\left[  g\right]  \right)  \leq Q\left(
\overline{\mathbb{B}^{n}},\mathbb{S}^{n-1}\right)  $ and moreover the infimum
is achieved if the inequality is strict. The strategy to solve the type II
Yamabe problem is to show that the strict inequality $Q\left(  M,\Sigma
,\left[  g\right]  \right)  <Q\left(  \overline{\mathbb{B}^{n}},\mathbb{S}%
^{n-1}\right)  $ is always true unless $\left(  M,\left[  g\right]  \right)  $
is conformal diffeomorphic to $\overline{\mathbb{B}^{n}}$. It has been
confirmed in various cases (see \cite{E2},\cite{M1} and \cite{M2}). But there
are still cases that remain open.

Apart from the minimization problem, both $Y\left(  M,\left[  g\right]
\right)  $ and $Q\left(  M,\Sigma,\left[  g\right]  \right)  $ are important
invariants and it is useful to have lower estimates for them. Let $\left(
X^{n},g_{+}\right)  $ be a Poincar\'{e}--Einstein manifold and $\Sigma
=\partial\overline{X}$. We pick a fixed defining function $r$ on $\overline
{X}$ which gives rise to a metric $\overline{g}=r^{2}g_{+}$ on $\overline{X}$.
As $\left[  \overline{g}\right]  $ and $\left[  \overline{g}|_{\Sigma}\right]
$ are invariantly defined, the Yamabe invariants $Y\left(  \overline
{X},\left[  \overline{g}\right]  \right)  ,Q\left(  \overline{X}%
,\Sigma,\left[  \overline{g}\right]  \right)  $ and $Y\left(  \Sigma,\left[
\overline{g}|_{\Sigma}\right]  \right)  $ are natural invariants of $\left(
X^{n},g_{+}\right)  $. X. Chen, M. Lai and F. Wang proved the following
elegant inequality relating these two Yamabe invariants.

\begin{theorem}
(Chen-Lai-Wang \cite{CLW}) Let $\left(  X^{n},g_{+}\right)  $ be a
Poincar\'{e}--Einstein manifold s.t. If the type II Yamabe problem on $\left(
\overline{X},\overline{g}\right)  $ has a minimizing solution, then
\begin{align*}
Y\left(  \Sigma,\left[  \overline{g}|_{\Sigma}\right]  \right)   &  \leq
\frac{n-2}{4\left(  n-1\right)  }Q\left(  \overline{X},\Sigma,\left[
\overline{g}\right]  \right)  ^{2},\text{ if }n\geq4;\\
32\pi\chi\left(  \Sigma\right)   &  \leq Q\left(  \overline{X},\Sigma,\left[
\overline{g}\right]  \right)  ^{2},\text{ if }n=3.
\end{align*}
Moreover, the equality holds if and only if $\left(  X^{n},g_{+}\right)  $ is
isometric to the hyperbolic space $\left(  \mathbb{H}^{n},g_{\mathbb{H}%
}\right)  $.
\end{theorem}

In our previous work \cite{WW}, we removed the restriction in Theorem 1 and
proved that the inequality is valid for all Poincar\'{e}--Einstein manifolds.
Since the inequality is vacuous when $Y\left(  \partial X,\left[  g\right]
\right)  \leq0$, we prefer to state the result in the following way.

\begin{theorem}
Let $\left(  X^{n},g_{+}\right)  $ be a Poincar\'{e}--Einstein manifold whose
conformal infinity has nonnegative Yamabe invariant. Then%
\begin{align*}
Q\left(  \overline{X},\Sigma,\left[  \overline{g}\right]  \right)   &
\geq2\sqrt{\frac{\left(  n-1\right)  }{\left(  n-2\right)  }Y\left(
\Sigma,\left[  \overline{g}|_{\Sigma}\right]  \right)  }\text{ if }n\geq4;\\
Q\left(  \overline{X},\Sigma,\left[  \overline{g}\right]  \right)   &
\geq4\sqrt{2\pi\chi\left(  \Sigma\right)  }\text{ if }n=3.
\end{align*}
Moreover, the equality holds iff $\left(  X^{n},g_{+}\right)  $ is isometric
to the hyperbolic space $\left(  \mathbb{H}^{n},g_{\mathbb{H}}\right)  $.
\end{theorem}

In this paper we prove that the same inequality holds in a much broader
context. It suffices for $\left(  X^{n},g_{+}\right)  $ to have Ricci
curvature bounded from below $Ric\left(  g_{+}\right)  \geq-\left(
n-1\right)  g_{+}$ and satisfy an asymptotic condition near infinity. This
seems to us to be the natural setting for the inequality and it fits well
within the general framework of understanding the boundary effect under a
Ricci curvature lower bound. We now explain the asymptotic condition
precisely. Let $\left(  X^{n},g_{+}\right)  $ be a conformally compact
manifold. As usual, we pick a fixed defining function $r$ on $\overline{X}$
which gives rise to a metric $\overline{g}=r^{2}g_{+}$ on $\overline{X}$. We
say that $\left(  X^{n},g_{+}\right)  $ is asymptotically Poincare-Einstein
if
\[
Ric\left(  g_{+}\right)  +\left(  n-1\right)  g_{+}=o\left(  r^{2}\right)  .
\]
We can now state our main result.

\begin{theorem}
\label{main}Let $\left(  X^{n},g_{+}\right)  $ be a conformally compact
manifold whose conformal infinity has nonnegative Yamabe invariant. If
$Ric\left(  g_{+}\right)  \geq-\left(  n-1\right)  g_{+}$ and $\left(
X^{n},g_{+}\right)  $ is asymptotically Poincare-Einstein, then
\begin{align*}
Q\left(  \overline{X},\Sigma,\left[  \overline{g}\right]  \right)   &
\geq2\sqrt{\frac{\left(  n-1\right)  }{\left(  n-2\right)  }Y\left(
\Sigma,\left[  \overline{g}|_{\Sigma}\right]  \right)  }\text{ if }n\geq4;\\
Q\left(  \overline{X},\Sigma,\left[  \overline{g}\right]  \right)   &
\geq4\sqrt{2\pi\chi\left(  \Sigma\right)  }\text{ if }n=3.
\end{align*}
Moreover, the equality holds iff $\left(  X^{n},g_{+}\right)  $ is isometric
to the hyperbolic space $\left(  \mathbb{H}^{n},g_{\mathbb{H}}\right)  $.
\end{theorem}

\begin{remark}
When $Y\left(  \Sigma,\left[  \overline{g}|_{\Sigma}\right]  \right)
=Y\left(  \mathbb{S}^{n-1}\right)  =\left(  n-1\right)  \left(  n-2\right)
\omega_{n-1}^{2/\left(  n-1\right)  }$, here $\omega_{n-1}$ is the volume of
$\mathbb{S}^{n-1}$, the right hand side then equals $2\left(  n-1\right)
\omega_{n-1}^{1/\left(  n-1\right)  }=Q\left(  \overline{\mathbb{B}^{n}%
},\mathbb{S}^{n-1}\right)  $. Thus in this case we must have equality and
hence rigidity. This rigidity results was proved by \cite{DJ} and \cite{LQS}.
Therefore our ineqaulity can be viewd as a quantative version of their
rigidity result: when the conformal infinity is closed to $\mathbb{S}^{n-1}$
in terms of the Yamabe invariant, $\left(  \overline{X},\left[  \overline
{g}\right]  \right)  $ is close to the ball $\overline{\mathbb{B}^{n}}$ in
terms of the type II Yamabe invariant.
\end{remark}

The method in \cite{CLW} is based on ideas introduced in Gursky-Han \cite{GH}
in which they studied the type I Yamabe invariant on $\overline{X}$. Let
$g\in\left[  \overline{g}\right]  $ be a type II Yamabe minimizer and write
$g_{+}=\rho^{-2}g$. The following identity plays an important role in the
proof of Theorem 1 as well as Theorem 2%
\[
T_{+}=T+\left(  n-2\right)  \rho^{-1}\left(  D^{2}\rho-\frac{\Delta\rho}%
{n}g\right)  ,
\]
where $T_{+}$ and $T$ are the traceless Ricci tensor of $g_{+}$ and $g$,
respectively. As $g_{+}$ is Einstein, $T_{+}=0$ and hence
\[
\frac{1}{n-2}\rho T=-\left(  D^{2}\rho-\frac{\Delta\rho}{n}g\right)  .
\]
By an integration by part over $X_{\varepsilon}=\left\{  r\geq\varepsilon
\right\}  $, using the fact that $g$ has constant scalar curvature, we obtain
\[
\frac{1}{n-2}\int_{X_{\varepsilon}}\rho\left\vert T\right\vert ^{2}%
dv_{g}=-\int_{\partial X_{\varepsilon}}T\left(  \nabla\rho,\nu\right)
d\sigma_{g}.
\]
The rest of the proof is to analyze the limit of the boundary term as
$\varepsilon\rightarrow0$.

When $g_{+}$ is not Einstein, the above approach breaks down at the beginning.
Instead, we study a modified Yamabe problem which produces a positive function
$u$ satisfying the equation
\[
-\Delta_{g_{+}}u=\frac{n\left(  n-2\right)  }{4}u.
\]
Write $u=v^{-\left(  n-2\right)  /2}$ and set $\Phi=v^{-1}\left(  \left\vert
\nabla v\right\vert ^{2}-v^{2}\right)  $. The following calculation is crucial
for our proof%
\[
\mathrm{div}\left(  v^{-\left(  n-2\right)  }\nabla\Phi\right)  =2v^{-\left(
n-2\right)  }Q,
\]
where
\[
Q=\left\vert D^{2}v-\frac{\Delta v}{n}g_{+}\right\vert ^{2}+Ric\left(  \nabla
v,\nabla v\right)  +\left(  n-1\right)  \left\vert \nabla v\right\vert ^{2}.
\]
We then integrate the above identity over $X_{\varepsilon}$. The analysis of
the boundary term follows the same strategy in \cite{CLW}.

The paper is organized as follows. In Section 2 we discuss some background
material. In Section 3, we study a modified Yamabe problem and estimate the
corresponding invariants. As a corollary we prove Theorem \ref{main}. We
discuss the related problem on compact manifolds in the last Section.

\section{Preliminaries}

Throughout this paper $\left(  X^{n},g_{+}\right)  $ is asymptotically
hyperbolic of order $C^{m,\alpha}$: if $r$ is smooth defining function on
$\overline{X}$, the metric $\overline{g}=r^{2}g$ extends to a $C^{m,\alpha}$
metric on $\overline{X}$ and $\left\vert d\rho\right\vert _{\overline{g}}%
^{2}=1$ along $\Sigma:=\partial\overline{X}$. For all the analysis it suffices
to have $m\geq4$. We also assume%
\[
Ric\left(  g_{+}\right)  \geq-\left(  n-1\right)  g_{+}%
\]
and that $g_{+}$ is asymptotically Poincare-Einstein in the following sense
\[
Ric\left(  g_{+}\right)  +\left(  n-1\right)  g_{+}=o\left(  r^{2}\right)  .
\]
Let $h\in\left[  \overline{g}|_{\Sigma}\right]  $ be a metric on $\Sigma$. It
is proved in \cite{Lee} that there is a defining function $r$ s.t. in a collar
neighborhood of $\Sigma$%
\begin{equation}
g_{+}=r^{-2}\left(  dr^{2}+h_{r}\right)  ,\label{gexp}%
\end{equation}
where $h_{r}$ is an $r$-dependent family of metrics on $\partial\overline{X}$
with $h_{r}|_{r=0}=h$. Moreover we have the following expansion (see, e.g.
\cite{GW})%
\[
h_{r}=h+h_{2}r^{2}+o\left(  r^{2}\right)  ,
\]
where%
\[
h_{2}=\left\{
\begin{array}
[c]{cc}%
-\frac{1}{n-3}\left(  Ric\left(  h\right)  -\frac{R_{h}}{2\left(  n-2\right)
}h\right)  , & \text{if }n\geq4;\\
-\frac{1}{4}h, & \text{if }n=3.
\end{array}
\right.
\]
It follows that $\overline{g}=r^{2}g_{+}$ has totally geodesic boundary. As we
assume $Y\left(  \Sigma,\left[  \overline{g}|_{\Sigma}\right]  \right)  \geq
0$, we choose $h$ to have scalar curvature $R_{h}\geq0$.

Lee \cite{Lee} constructed a positive smooth function $\phi$ on $X$ s.t.
$\Delta\phi=n\phi$ and near $\partial\overline{X}$%
\[
\phi=r^{-1}+\frac{R_{h}}{4\left(  n-1\right)  \left(  n-2\right)  }r+o\left(
r^{2}\right)  .
\]
Under the condition $R_{h}\geq0$, he further proved that $\left\vert
d\phi\right\vert _{g_{+}}^{2}\leq\phi^{2}$. Consider the metric $\widetilde
{g}:=\phi^{-2}g_{+}$ on $\overline{X}$ . Its scalar curvature is given by%
\begin{align}
\widetilde{R}  &  =\phi^{2}\left(  R+2\left(  n-1\right)  \phi^{-1}\Delta
\phi-n\left(  n-1\right)  \phi^{-2}\left\vert d\phi\right\vert ^{2}\right)
\label{sali}\\
&  \geq\phi^{2}\left(  R+n\left(  n-1\right)  \right)  .\nonumber
\end{align}
Moreover, by a direct calculation the boundary is totally geodesic. We
consider the following modified energy functional%
\[
\widetilde{E}\left(  f\right)  =E_{\overline{g}}\left(  f\right)  -\int
_{X}\left(  R_{+}+n\left(  n-1\right)  \right)  \phi^{2}f^{2}dv_{\overline{g}%
}.
\]
Note that $\left(  R+n\left(  n-1\right)  \right)  \phi^{2}\in C^{m-3,\alpha
}\left(  \overline{X}\right)  $ under our assumptions. More explicitly, by
(\ref{sali})
\[
\widetilde{E}\left(  f\right)  =\int_{X}\left[  \frac{4\left(  n-1\right)
}{n-2}\left\vert df\right\vert _{\widetilde{g}}^{2}+\left(  \widetilde
{R}-\left(  R+n\left(  n-1\right)  \right)  \phi^{2}\right)  f^{2}\right]
dv_{\overline{g}}\geq0.
\]
Since $R_{+}+n\left(  n-1\right)  \geq0$, we have
\begin{equation}
E_{\overline{g}}\left(  f\right)  \geq\widetilde{E}\left(  f\right)  .
\label{c2e}%
\end{equation}

\section{Estimate on modified Yamabe Quotients}

For $1<q\leq n/\left(  n-2\right)  $, consider%
\[
\widetilde{\lambda}_{q}:=\inf\frac{\widetilde{E}\left(  f\right)  }{\left(
\int_{\Sigma}\left\vert f\right\vert ^{q+1}d\sigma_{\overline{g}}\right)
^{2/\left(  q+1\right)  }}.
\]

\begin{theorem}
\label{isub}Let $\left(  X^{n},g_{+}\right)  $ be a Poincar\'{e}--Einstein
manifold whose conformal infinity has positive Yamabe invariant. For $1<q\leq
n/\left(  n-2\right)  $ the invariant $\lambda_{q}$ satisfies%
\begin{align*}
\widetilde{\lambda}_{q}  &  \geq2\sqrt{\frac{\left(  n-1\right)  }{\left(
n-2\right)  }Y\left(  \Sigma,\left[  \overline{g}|_{\Sigma}\right]  \right)
}V\left(  \Sigma,\overline{g}\right)  ^{-\frac{\left(  n-q\left(  n-2\right)
\right)  }{\left(  n-3\right)  \left(  q+1\right)  }}\text{ if }n\geq4;\\
\widetilde{\lambda}_{q}  &  \geq4\sqrt{2\pi\chi\left(  \Sigma\right)
}V\left(  \Sigma,\overline{g}\right)  ^{-\frac{3-q}{2\left(  q+1\right)  }%
}\text{ if }n=3.
\end{align*}

\end{theorem}

Since $\widetilde{E}\left(  f\right)  \geq0$, it is easy to see that
$\lim_{q\nearrow n/\left(  n-2\right)  }\widetilde{\lambda}_{q}=\widetilde
{\lambda}_{n/\left(  n-2\right)  }$. Therefore, it suffices to prove the above
theorem for $q<n/\left(  n-2\right)  $.

Since the trace operator $H^{1}\left(  \overline{X}\right)  \rightarrow
L^{q+1}\left(  \Sigma\right)  $ is compact for $q<n/\left(  n-2\right)  $, by
standard elliptic theory, the above infimum $\lambda_{q}$ is achieved by a
smooth, positive function $f$ s.t.%
\begin{equation}
\int_{\Sigma}f^{q+1}d\overline{\sigma}=1 \label{Norm}%
\end{equation}
and%
\begin{equation}
\left\{
\begin{array}
[c]{cc}%
-\frac{4\left(  n-1\right)  }{n-2}\overline{\Delta}f+\overline{R}f=\left(
R+n\left(  n-1\right)  \right)  \phi^{2}f & \text{on }\overline{X},\\
\frac{4\left(  n-1\right)  }{n-2}\frac{\partial f}{\partial\overline{\nu}%
}=\lambda_{q}f^{q} & \text{on }\Sigma.
\end{array}
\right.  \label{keq}%
\end{equation}
By the conformal invariance of the conformal Laplacian, we have%
\begin{align*}
L_{g}\left(  f\phi^{-\left(  n-2\right)  /2}\right)   &  =\phi^{-\left(
n+2\right)  /2}L_{\overline{g}}\left(  f\right) \\
&  =\left(  R+n\left(  n-1\right)  \right)  f\phi^{-\left(  n-2\right)  /2}.
\end{align*}
In other words, $u:=f\phi^{-\left(  n-2\right)  /2}$ satisfies the following
equation
\begin{equation}
-\Delta_{g_{+}}u=\frac{n\left(  n-2\right)  }{4}u. \label{efu}%
\end{equation}
Write $u=v^{-\left(  n-2\right)  /2}$. Then%
\[
\Delta_{g_{+}}v=\frac{n}{2}v^{-1}\left(  \left\vert dv\right\vert _{g_{+}}%
^{2}+v^{2}\right)  .
\]
Equivalently $\Delta_{g_{+}}v-nv=\frac{n}{2}\Phi$ with $\Phi=v^{-1}\left(
\left\vert dv\right\vert _{g_{+}}^{2}-v^{2}\right)  $.

\begin{lemma}
We have
\begin{equation}
\mathrm{div}\left(  v^{-\left(  n-2\right)  }\nabla\Phi\right)  =2v^{-\left(
n-2\right)  }Q, \label{ki}%
\end{equation}
where
\[
Q=\left\vert D^{2}v-\frac{\Delta v}{n}g_{+}\right\vert ^{2}+Ric\left(  \nabla
v,\nabla v\right)  +\left(  n-1\right)  \left\vert \nabla v\right\vert
^{2}\geq0.
\]
All the computation is done with respect to $g_{+}$, but we drop the subscript
to simplify the presentation.
\end{lemma}

\begin{proof}
As $v\Phi=\left\vert \nabla v\right\vert ^{2}-v^{2}$, we have, by using the
Bochner formula
\begin{align*}
\frac{1}{2}\left(  v\Delta\phi+2\left\langle \nabla v,\nabla\phi\right\rangle
+\phi\Delta v\right)   &  =\left\vert D^{2}v\right\vert ^{2}+\left\langle
\nabla v,\nabla\Delta v\right\rangle +Ric\left(  \nabla v,\nabla v\right)
-v\Delta v-\left\vert \nabla v\right\vert ^{2}\\
&  =\frac{\left(  \Delta v\right)  ^{2}}{n}+\left\langle \nabla v,\nabla\Delta
v\right\rangle +v\Delta v-n\left\vert \nabla v\right\vert ^{2}+Q\\
&  =\frac{\Delta v}{n}\left(  \Delta v-nv\right)  +\left\langle \nabla
v,\nabla\left(  \Delta v-nv\right)  \right\rangle +Q\\
&  =\frac{1}{2}\Phi\Delta v+\frac{n}{2}\left\langle \nabla v,\nabla
\Phi\right\rangle +Q
\end{align*}
Thus,%
\[
\Delta\Phi=\left(  n-2\right)  v^{-1}\left\langle \nabla v,\nabla
\Phi\right\rangle +2Q
\]
or
\[
\mathrm{div}\left(  v^{-\left(  n-2\right)  }\nabla\Phi\right)  =2v^{-\left(
n-2\right)  }Q.
\]

\end{proof}

We now consider the metric $g=u^{4/\left(  n-2\right)  }g_{+}$. Since
$u=f\phi^{-\left(  n-2\right)  /2}$, we also have%
\[
g=f^{4/\left(  n-2\right)  }\phi^{-2}g_{+}=f^{4/\left(  n-2\right)
}\widetilde{g}.
\]
As $\partial\overline{X}$ is totally geodesic w.r.t. $\widetilde{g}$ and $g$
is conformal to $\widetilde{g}$, we know that $\partial\overline{X}$ is
umbilic w.r.t. $g$ and its mean curvature, in view of the boundary condition
of (\ref{keq}), is given by
\begin{equation}
H=\frac{\lambda_{q}}{2}f^{q-\frac{n}{n-2}}. \label{Hf}%
\end{equation}
Set $\rho=u^{2/\left(  n-2\right)  }=v^{-1}$. By a direct calculation, the
equation (\ref{efu}) becomes, using $g$ as the background metric
\begin{equation}
2\rho\Delta\rho=n\left(  \left\vert \nabla\rho\right\vert ^{2}-1\right)  .
\label{sce}%
\end{equation}
Let $t$ be the geodesic distance to $\Sigma$ w.r.t. $g$. We need the following
lemma which is essentially contained in \cite{CLW}.

\begin{lemma}
\label{bda}\bigskip Near $\Sigma=\partial\overline{X}$, we can write
\[
g=dt^{2}+g_{ij}\left(  t,x\right)  dx_{i}dx_{j},
\]
where $\left\{  x_{1},\cdots,x_{n-1}\right\}  $ are local coordinates on
$\Sigma$. Then
\[
\rho=t-\frac{H}{2\left(  n-1\right)  }t^{2}+\frac{1}{6}\left(  \frac
{R^{\Sigma}}{n-2}-\frac{H^{2}}{n-1}\right)  t^{3}+o\left(  t^{3}\right)  .
\]
In particular,%
\[
\frac{\partial}{\partial\nu}\left[  \rho^{-1}\left(  \left\vert \nabla
\rho\right\vert ^{2}-1\right)  \right]  |_{\Sigma}=\frac{R^{\Sigma}}%
{n-2}-\frac{H^{2}}{n-1}.
\]

\end{lemma}

\begin{proof}
For completeness, we present the proof showing that the Einstein condition is
not required. In local coordinates%
\begin{align*}
\left\vert \nabla\rho\right\vert ^{2}  &  =\left(  \frac{\partial\rho
}{\partial t}\right)  ^{2}+g^{ij}\frac{\partial\rho}{\partial x_{i}}%
\frac{\partial\rho}{\partial x_{j}},\\
\Delta\rho &  =\frac{\partial^{2}\rho}{\partial t^{2}}+\frac{\partial\log
\sqrt{G}}{\partial t}\frac{\partial\rho}{\partial t}+\frac{1}{\sqrt{G}}%
\frac{\partial}{\partial x_{i}}\left(  g^{ij}\sqrt{G}\frac{\partial\rho
}{\partial x_{j}}\right)  .
\end{align*}
Restricting (\ref{sce}) on $\Sigma$ on which both $\rho$ and $r$ vanish with
order $1$ yields $\frac{\partial\rho}{\partial t}|_{\Sigma}=1$.

Differentiating (\ref{sce}) in $t$ yields
\begin{equation}
\frac{2}{n}\left(  \frac{\partial\rho}{\partial t}\Delta\rho+\rho
\frac{\partial}{\partial t}\Delta\rho\right)  =2\frac{\partial\rho}{\partial
t}\frac{\partial^{2}\rho}{\partial t^{2}}+2g^{ij}\frac{\partial^{2}\rho
}{\partial x_{i}\partial t}\frac{\partial\rho}{\partial x_{j}}-g^{ik}%
g^{jl}\frac{\partial g_{kl}}{\partial t}\frac{\partial\rho}{\partial x_{i}%
}\frac{\partial\rho}{\partial x_{j}}. \label{d1}%
\end{equation}
Evaluating both sides on $\Sigma$ yields%
\[
\frac{2}{n}\left(  \frac{\partial^{2}\rho}{\partial t^{2}}+\frac{\partial
\log\sqrt{G}}{\partial t}\right)  |_{\Sigma}=2\frac{\partial^{2}\rho}{\partial
t^{2}}|_{\Sigma}.
\]
Thus%
\[
\frac{\partial^{2}\rho}{\partial t^{2}}|_{\Sigma}=\frac{1}{n-1}\frac
{\partial\log\sqrt{G}}{\partial t}|_{\Sigma}=-\frac{H}{n-1}.
\]
Differentiating the formula for $\Delta\rho$ we get%
\begin{align*}
\frac{\partial}{\partial t}\Delta\rho|_{\Sigma}  &  =\left(  \frac
{\partial^{3}\rho}{\partial t^{3}}+\frac{\partial^{2}\log\sqrt{G}}{\partial
t^{2}}+\frac{\partial\log\sqrt{G}}{\partial t}\frac{\partial^{2}\rho}{\partial
t^{2}}\right)  |_{\Sigma}\\
&  =\left(  \frac{\partial^{3}\rho}{\partial t^{3}}+\frac{\partial^{2}%
\log\sqrt{G}}{\partial t^{2}}+\frac{H^{2}}{n-1}\right)  |_{\Sigma}%
\end{align*}
Differentiating (\ref{d1}) in $r$ and evaluating on $\Sigma$, we obtain%
\[
\frac{2}{n}\left(  \frac{\partial^{2}\rho}{\partial t^{2}}\Delta\rho
+2\frac{\partial}{\partial t}\Delta\rho\right)  |_{\Sigma}=2\left(
\frac{\partial^{2}\rho}{\partial t^{2}}\right)  ^{2}|_{\Sigma}+2\frac
{\partial^{3}\rho}{\partial t^{3}}|_{\Sigma}=\frac{2H^{2}}{\left(  n-1\right)
^{2}}+2\frac{\partial^{3}\rho}{\partial t^{3}}|_{\Sigma}.
\]
Using the previous formulas, we arrive at
\[
\frac{\partial^{3}\rho}{\partial t^{3}}|_{\Sigma}=\frac{2}{n-2}\left(
\frac{H^{2}}{n-1}+\frac{\partial^{2}\log\sqrt{G}}{\partial t^{2}}|_{\Sigma
}\right)  .
\]
By a direct calculation, we also have%
\[
\frac{\partial^{2}\log\sqrt{G}}{\partial t^{2}}|_{\Sigma}=-Ric\left(  \nu
,\nu\right)  -\frac{H^{2}}{n-1}.
\]
Therefore%
\begin{align*}
\frac{\partial^{3}\rho}{\partial t^{3}}|_{\Sigma}  &  =-\frac{2}%
{n-2}Ric\left(  \nu,\nu\right) \\
&  =\frac{R^{\Sigma}}{n-2}-\frac{H^{2}}{n-1},
\end{align*}
where we used the Gauss equation in the last step.

The second identity follows from a direct calculation.
\end{proof}

\bigskip

We can now prove Theorem \ref{isub}. Integrating the identity (\ref{ki}) on
$X_{\varepsilon}=\left\{  t\geq\varepsilon\right\}  $ yields%
\[
2\int_{X_{\varepsilon}}v^{-\left(  n-2\right)  }Qdv_{g_{+}}=\int_{\partial
X_{\varepsilon}}v^{-\left(  n-2\right)  }\frac{\partial\Phi}{\partial\nu
}d\sigma_{g_{+}}.
\]
Since $g_{+}=\rho^{-2}g$, we obtain by a direct calculation%
\[
\int_{\partial X_{\varepsilon}}v^{-\left(  n-2\right)  }\frac{\partial\Phi
}{\partial\nu^{+}}d\sigma_{g_{+}}=\int_{\partial X_{\varepsilon}}%
\frac{\partial}{\partial\nu}\left[  \rho^{-1}\left(  \left\vert \nabla
\rho\right\vert ^{2}-1\right)  \right]  d\sigma_{g}.
\]
Therefore%
\[
2\int_{X_{\varepsilon}}v^{-\left(  n-2\right)  }Qdv_{g_{+}}=\int_{\partial
X_{\varepsilon}}\frac{\partial}{\partial\nu}\left[  \rho^{-1}\left(
\left\vert \nabla\rho\right\vert ^{2}-1\right)  \right]  d\sigma_{g}.
\]
Letting $\varepsilon\rightarrow0$, we obtain, in view of Lemma \ref{bda}%
\begin{equation}
2\int_{X}v^{-\left(  n-2\right)  }Qdv_{g_{+}}=\int_{\Sigma}\left(
\frac{R^{\Sigma}}{n-2}-\frac{H^{2}}{n-1}\right)  d\sigma_{g} \label{fe1}%
\end{equation}

\ \ The rest of the argument is the same as in \cite{WW}. We present it for
completeness. By (\ref{Hf}) and the Holder inequality again%
\begin{align*}
\int_{\Sigma}H^{2}d\sigma &  =\left(  \frac{\lambda_{q}}{2}\right)  ^{2}%
\int_{\Sigma}f^{2\left(  q-\frac{n}{n-2}\right)  }f^{2\left(  n-1\right)
/\left(  n-2\right)  }d\overline{\sigma}\\
&  =\left(  \frac{\lambda_{q}}{2}\right)  ^{2}\int_{\Sigma}f^{2\left(
q-\frac{1}{n-2}\right)  }d\overline{\sigma}\\
&  \leq\left(  \frac{\lambda_{q}}{2}\right)  ^{2}\left(  \int_{\Sigma}%
f^{q+1}d\overline{\sigma}\right)  ^{2\left(  q-\frac{1}{n-2}\right)  /\left(
q+1\right)  }V\left(  \Sigma,\overline{g}\right)  ^{\left(  \frac{n}%
{n-2}-q\right)  /\left(  q+1\right)  }\\
&  =\left(  \frac{\lambda_{q}}{2}\right)  ^{2}V\left(  \Sigma,\overline
{g}\right)  ^{\left(  \frac{n}{n-2}-q\right)  /\left(  q+1\right)  }.
\end{align*}
Plugging the above inequality into (\ref{fe1}), we obtain%
\begin{equation}
2\int_{X_{\varepsilon}}v^{-\left(  n-2\right)  }Qdv_{g_{+}}\leq\frac
{\lambda_{q}^{2}}{4\left(  n-1\right)  }V\left(  \Sigma,\overline{g}\right)
^{\left(  \frac{n}{n-2}-q\right)  /\left(  q+1\right)  }-\frac{1}{n-2}%
\int_{\Sigma}R^{\Sigma}d\sigma. \label{fe2}%
\end{equation}
When $n=3$, this implies%
\[
\lambda_{q}^{2}V\left(  \Sigma,\overline{g}\right)  ^{\left(  3-q\right)
/\left(  q+1\right)  }\geq32\pi\chi\left(  \Sigma\right)  .
\]
In the following, we assume $n>3\,$. By (\ref{Norm}) and the Holder
inequality
\begin{align*}
1  &  =\int_{\Sigma}f^{q+1}d\overline{\sigma}\\
&  \leq\left(  \int_{\Sigma}f^{2\left(  n-1\right)  /\left(  n-2\right)
}d\overline{\sigma}\right)  ^{\frac{\left(  q+1\right)  \left(  n-2\right)
}{2\left(  n-1\right)  }}V\left(  \Sigma,\overline{g}\right)  ^{\frac
{n-q\left(  n-2\right)  }{2\left(  n-1\right)  }}\\
&  =V\left(  \Sigma,g\right)  ^{\frac{\left(  q+1\right)  \left(  n-2\right)
}{2\left(  n-1\right)  }}V\left(  \Sigma,\overline{g}\right)  ^{\frac
{n-q\left(  n-2\right)  }{2\left(  n-1\right)  }}%
\end{align*}
Thus%
\[
V\left(  \Sigma,\overline{g}\right)  ^{-\frac{n-q\left(  n-2\right)  }{\left(
n-2\right)  \left(  q+1\right)  }}\leq V\left(  \Sigma,g\right)  .
\]
Plugging this inequality into (\ref{fe2}) yields
\begin{align*}
&  2\int_{X_{\varepsilon}}v^{-\left(  n-2\right)  }Qdv_{g_{+}}\\
\leq &  \frac{V\left(  \Sigma,g\right)  ^{\frac{n-1}{n-3}}}{4\left(
n-1\right)  }\left[  \widetilde{\lambda}_{q}^{2}V\left(  \Sigma,\overline
{g}\right)  ^{\frac{2\left(  n-q\left(  n-2\right)  \right)  }{\left(
n-3\right)  \left(  q+1\right)  }}-\frac{4\left(  n-1\right)  }{\left(
n-2\right)  V\left(  \Sigma,g\right)  ^{\frac{n-1}{n-3}}}\int_{\Sigma
}R^{\Sigma}d\sigma\right] \\
\leq &  \frac{V\left(  \Sigma,g\right)  ^{\frac{n-1}{n-3}}}{4\left(
n-1\right)  }\left[  \widetilde{\lambda}_{q}^{2}V\left(  \Sigma,\overline
{g}\right)  ^{\frac{2\left(  n-q\left(  n-2\right)  \right)  }{\left(
n-3\right)  \left(  q+1\right)  }}-\frac{4\left(  n-1\right)  }{\left(
n-2\right)  }Y\left(  \Sigma,\left[  \gamma\right]  \right)  \right]  .
\end{align*}
Therefore
\[
\widetilde{\lambda}_{q}^{2}\geq\frac{4\left(  n-1\right)  }{\left(
n-2\right)  }Y\left(  \Sigma\right)  V\left(  \Sigma,\overline{g}\right)
^{-\frac{2\left(  n-q\left(  n-2\right)  \right)  }{\left(  n-3\right)
\left(  q+1\right)  }}.
\]
This finishes the proof of Theorem \ref{isub}

We are now ready to prove our main result.

\begin{theorem}
Let $\left(  X^{n},g_{+}\right)  $ be a conformally compact manifold whose
conformal infinity has nonnegative Yamabe invariant. If $Ric\left(
g_{+}\right)  \geq-\left(  n-1\right)  g_{+}$ and $\left(  X^{n},g_{+}\right)
$ is asymptotically Poincare-Einstein, then
\begin{align*}
Q\left(  \overline{X},\Sigma,\left[  \overline{g}\right]  \right)   &
\geq2\sqrt{\frac{\left(  n-1\right)  }{\left(  n-2\right)  }Y\left(
\Sigma,\left[  \overline{g}|_{\Sigma}\right]  \right)  }\text{ if }n\geq4;\\
Q\left(  \overline{X},\Sigma,\left[  \overline{g}\right]  \right)   &
\geq4\sqrt{2\pi\chi\left(  \Sigma\right)  }\text{ if }n=3.
\end{align*}
Moreover, the equality holds iff $\left(  X^{n},g_{+}\right)  $ is isometric
to the hyperbolic space $\left(  \mathbb{H}^{n},g_{\mathbb{H}}\right)  $.
\end{theorem}

\begin{proof}
By (\ref{c2e}), we have $Q\left(  \overline{X},\Sigma,\left[  \overline
{g}\right]  \right)  \geq\widetilde{\lambda}_{n/\left(  n-2\right)  }$.
Therefore the inequality follows immediately from Theorem \ref{isub}.

Suppose the equality holds. We present the argument for $n\geq4$ and the same
argument works for $n=3$ with trivial modification. If $Y\left(
\Sigma,\left[  \overline{g}|_{\Sigma}\right]  \right)  <Y\left(
\mathbb{S}^{n-1}\right)  $, the equality then implies%
\[
Q\left(  \overline{X},\Sigma,\left[  \overline{g}\right]  \right)
=\widetilde{\lambda}_{n/\left(  n-2\right)  }<Q\left(  \overline
{\mathbb{B}^{n}},\mathbb{S}^{n-1}\right)  .
\]
Just like in the original Yamabe problem, this strict inequality implies that
$\widetilde{\lambda}_{n/\left(  n-2\right)  }$ is achieved. Therefore in the
proof of Theorem \ref{isub}, we can take $q=n/\left(  n-2\right)  $ and obtain%
\[
2\int_{X_{\varepsilon}}v^{-\left(  n-2\right)  }Qdv_{g_{+}}\leq\frac{V\left(
\Sigma,g\right)  ^{\frac{n-1}{n-3}}}{4\left(  n-1\right)  }\left[
\widetilde{\lambda}_{n/\left(  n-2\right)  }^{2}-\frac{4\left(  n-1\right)
}{\left(  n-2\right)  }Y\left(  \Sigma,\left[  \gamma\right]  \right)
\right]  =0.
\]
Thus $Q=0$. In particular, $v>0$ satisfies the over-determined system%
\[
D^{2}v=\frac{\Delta v}{n}g_{+}.
\]
This implies that $\left(  X^{n},g_{+}\right)  $ is isometric to the
hyperbolic space (cf. \cite{CLW} for the argument).

If $Y\left(  \Sigma,\left[  \overline{g}|_{\Sigma}\right]  \right)  =Y\left(
\mathbb{S}^{n-1}\right)  $, then $\left(  \Sigma,\left[  \overline{g}%
|_{\Sigma}\right]  \right)  $ is conformally equivalent to $\mathbb{S}^{n-1}$,
by the solution of the Yamabe problem for closed manifolds. Then $\left(
X^{n},g_{+}\right)  $ is isometric to the hyperbolic space $\left(
\mathbb{H}^{n},g_{\mathbb{H}}\right)  $ by \cite{DJ} and \cite{LQS}.
\end{proof}

\section{Some Discussions on Compact Manifolds with Boundary}

It is a natural question if the inequality holds for a compact Riemannian
manifold $\left(  M^{n},g\right)  $ with $Ric$ and $\Pi\geq1$. We are
motivated by the observation that some results for conformally compact
manifolds follow from results for compact Riemannian manifolds by a limiting
process. As an illustration, consider the following theorem by Lee.

\begin{theorem}
(Lee \cite{Lee}) Let $\left(  X^{n},g_{+}\right)  $ be a conformally compact
manifold whose conformal infinity has nonnegative Yamabe invariant. If
$Ric\left(  g_{+}\right)  \geq-\left(  n-1\right)  g_{+}$ and $\left(
X^{n},g_{+}\right)  $ is asymptotically Poincare-Einstein, then the bottom of
spectrum $\lambda_{0}\left(  X^{n},g_{+}\right)  =\left(  n-1\right)  ^{2}/4$.
\end{theorem}

When the Yambabe invariant of the conformal infinity is positive, Lee's
theorem follows from the following result for compact Riemannian manifolds.

\begin{theorem}
Let $\left(  M^{n},g\right)  $ be a compact Riemannian manifold with
$Ric\geq-\left(  n-1\right)  $. If along the boundary $\Sigma:=\partial M$ we
have the mean curvature $H\geq n-1$, then the first Dirichlet eigenvalue
\[
\lambda_{0}\left(  M\right)  \geq\frac{\left(  n-1\right)  ^{2}}{4}.
\]

\end{theorem}

Let $r$ be the distance function to $\Sigma$. By standard method in Riemannian
geometry, we have
\[
\Delta r\leq-\left(  n-1\right)
\]
in the support sense. A direct calculation yields
\[
\Delta e^{\left(  n-1\right)  r/2}\leq-\frac{\left(  n-1\right)  ^{2}}%
{4}e^{\left(  n-1\right)  r/2}.
\]
This implies $\lambda_{0}\left(  M\right)  \geq\frac{\left(  n-1\right)  ^{2}%
}{4}$ (for technical details see \cite{Wa1}).

We can deduce Lee's theorem from Theorem when the conformal infinity has
positive Yamabe invariant in the following way. As explained in Section 2, we
pick a metric $h$ on the conformal infinity with positive scalar curvature and
then we have a good defining function $r$ s.t. near the conformal infinity
$g_{+}$ has a nice expansion (\ref{gexp}). Then a simple calculation shows
that the mean curvature of the boundary of $X_{\varepsilon}:=\left\{
r\geq\varepsilon\right\}  $ satisfies%
\[
H=n-1+\frac{R_{h}}{2\left(  n-2\right)  }\varepsilon^{2}+o\left(
\varepsilon^{2}\right)  .
\]
As $R_{h}>0$, we have $H>n-1$ if $\varepsilon$ is small enough. By Theorem,
$\lambda_{0}\left(  X_{\varepsilon}\right)  \geq\frac{\left(  n-1\right)
^{2}}{4}$. It follows that $\lambda_{0}\left(  X\right)  \geq\frac{\left(
n-1\right)  ^{2}}{4}$. As the opposite inequality was known by \cite{Ma}, we
have $\lambda_{0}\left(  X\right)  =\frac{\left(  n-1\right)  ^{2}}{4}$. When
the conformal infinity has zero Yamabe invariant, the situation is more
subtle. But by an idea in Cai-Galloway\cite{CG}, a similar argument still
works (cf. \cite{Wa1}).

We now come back to Theorem. By the asymptotic expansion (\ref{gexp}) the
second fundamental form of $\partial X_{\varepsilon}$ satisfies%
\[
\Pi_{+}=\left(  1+O\left(  \varepsilon\right)  \right)  g_{+},
\]
i.e. all the principal curvatures are close to $1$. This leads us to consider
a compact Riemannian manifold $\left(  M^{n},g\right)  $ with $Ric\geq-\left(
n-1\right)  $ and $\Pi\geq1$ on its boundary $\Sigma$ and ask the question
whether the inequality
\begin{align}
Q\left(  M,\Sigma,g\right)   &  \geq2\sqrt{\frac{\left(  n-1\right)  }{\left(
n-2\right)  }Y\left(  \Sigma\right)  }\text{ if }n\geq4;\label{cje}\\
Q\left(  M,\Sigma,g\right)   &  \geq4\sqrt{2\pi\chi\left(  \Sigma\right)
}\text{ if }n=3\nonumber
\end{align}
holds. The answer turns out to be no in general. To construct counter example,
we consider the hyperbolic space using the ball model $\mathbb{B}^{n}$ with
the metric $g_{\mathbb{H}}=\frac{4}{\left(  1-\left\vert x\right\vert
^{2}\right)  ^{2}}dx^{2}$. For $0<R<1$, the Euclidean ball
\[
\left\{  x\in\mathbb{B}^{n}:\left\vert x\right\vert ^{2}=\sum_{i=1}^{n}%
x_{i}^{2}\leq R\right\}
\]
is a geodesic ball in $\left(  \mathbb{B}^{n},g_{\mathbb{H}}\right)  $ and the
boundary has 2nd fundamental form $\Pi=\frac{1+R^{2}}{2R}I$. We now consider
\[
M=\left\{  x\in\mathbb{B}^{n}:\left\vert x\right\vert ^{2}=\sum_{i=1}%
^{n-1}x_{i}^{2}+kx_{n}^{2}\leq R\right\}  ,
\]
where $k>0$ is close to $1$. Then $\left(  M,g_{\mathbb{H}}\right)  $ is a
compact hyperbolic manifold with boundary and on its boundary we have $\Pi
\geq1$ if $k$ is sufficiently close to $1$ by continuity. Since $\Sigma$ with
the induced metric is rotationally symmetric, it is conformally equivalent to
the standard sphere $\mathbb{S}^{n-1}$. Thus, $Y\left(  \Sigma\right)
=Y\left(  \mathbb{S}^{n-1}\right)  $. But when $k\neq1$, the boundary is not
umbilic with respect to the Euclidean metric and hence not with respect to
$g_{\mathbb{H}}$ either. By \cite{E2} and \cite{M2}, $Q\left(  M,\Sigma
,g_{\mathbb{H}}\right)  <Q\left(  \overline{\mathbb{B}^{n}},\mathbb{S}%
^{n-1}\right)  $. It follows that the inequality (\ref{cje}) is false.

Therefore, for a compact Riemannian manifold $\left(  M^{n},g\right)  $ with
$Ric\geq-\left(  n-1\right)  $ and $\Pi\geq1$ on its boundary $\Sigma$, it is
more subtle to estimate its type II Yamabe invariant in terms of the boundary
geometry. It is an interesting question and we do not have an explicit
conjecture. Let us mention that in a similar setting, namely for a compact
$\left(  M^{n},g\right)  $ with $Ric\geq0$ and $\Pi\geq1$ on its boundary
$\Sigma$, there is a well-formulated conjecture \cite{Wa2} on the type II
Yamabe invariant in terms of the boundary area.

\end{document}